\documentclass[%
 reprint,
 amsmath,amssymb,showkeys,
 aps,
]{revtex4-2}

\usepackage{graphicx,xcolor}
\usepackage{dcolumn}
\usepackage{bm}
\usepackage{hyperref}
\usepackage[mathlines]{lineno}

\begin{document}

\preprint{APS/123-QED}

\title{Power Law Behavior of Center-Like Decaying Oscillation : Exponent through Perturbation Theory and Optimization}

\author{Sandip Saha}

\email{sahasandip.loknath@gmail.com}
\email{sandip.saha@kaust.edu.sa}
\affiliation{King Abdullah University of Science and Technology, Thuwal 23955-6900, Saudi Arabia}%

\begin{abstract}
In dynamical systems theory, there is a lack of a straightforward rule to distinguish exact center solutions from decaying center-like solutions, as both require the damping force function to be zero~\cite{len3.5,powerlaw}. {\color{black}By adopting a multi-scale perturbative method, we have demonstrated a general rule for the decaying center-like power law behavior, characterized by an exponent of $\frac{1}{3}$. The investigation began with a physical question about the higher-order nonlinearity in a damping force function, which exhibits birhythmic and trirhythmic behavior under a transition to a decaying center-type solution. Using numerical optimization algorithms, we identified the power law exponent for decaying center-type behavior across various rhythmic conditions. For all scenarios, we consistently observed a decaying power law with an exponent of $\frac{1}{3}$.}Our study aims to elucidate their dynamical differences, contributing to theoretical insights and practical applications where distinguishing between different types of center-like behaviour is crucial. This key result would be beneficial for studying the multi-rhythmic nature of biological and engineering systems.
\end{abstract}

\maketitle


\section{Introduction}
Science is a common approach used to deepen our understanding of nature, encompassing core branches such as physics, chemistry, and biology.  Over time, new interdisciplinary fields like biophysics and biomedical have emerged, merging fundamental sciences with mathematical tools to validate findings through intuitive mathematical frameworks rather than solely relying on laboratory experiments or computer simulations~\cite{JeffHasty,2018_RNA_Velocity_Nature,2021_RNA_Velocity_MBS}.  Optimization techniques have integrated with these core sciences, often used for machine learning methods. It play a crucial role in studying both real and prototype systems through dynamic systems modeling, offering insights into system dynamics in real-world scenarios. When studying dynamic systems, oscillations frequently arise, characterized by isolated trajectories known as limit cycle oscillations, which demonstrate self-sustained periodic behavior~\cite{bender,murraynld,murraybio,sabatini}.  Examples include chemical oscillations~\cite{nicolis1977,epstein,grayscott,kuramoto} and biological systems that includes glycolytic oscillations~\cite{glyscott,goldbook,gly2,gly3}, calcium oscillations~\cite{cal}, cell division~\cite{chen2004}, and circadian rhythms~\cite{cirgoldbeter1995,sen2008}.Extensive literature exists on limit cycle oscillations, both theoretically and experimentally, illustrating their role as a well-understood natural phenomenon sustained by energy balance mechanisms~\cite{SahaCNSNS,SahaCNSNS2}.  For strongly nonlinear systems, analytical study methods are limited, often necessitating reliance on computer simulations.  In contrast, weakly nonlinear systems benefit from multi-scale perturbative methods like Renormalization Group (RG)~\cite{chen1,chen2,goto2007renormalization} and Krylov-Boguliuvov (KB)~\cite{slross,len0,SahaCNSNS,lccounting}, which effectively predict system behaviors, including analytical solutions~\cite{powerlaw} through amplitude and phase dynamics.  These perturbative methods provide insights into coupled and uncoupled amplitude-phase dynamics, crucial for understanding system behaviors.  A common technique involves introducing an artificial nonlinear parameter to study strongly nonlinear systems via perturbative methods, subsequently setting it to unity once amplitude and phase equations are derived.

Dynamical systems capable of exhibiting isochronous oscillations are crucial for modeling real-world systems. Isochronicity, which refers to oscillations with a constant period regardless of their amplitude, is not only significant for stability and bifurcation theory but also plays a key role in understanding limit cycles and oscillations with a center. Recent research has explored scenarios where a self-sustained limit cycle oscillation transitions suddenly into an isochronous state due to external nuisance of environment~\cite{len3.5,powerlaw}. In a study (cf.~\cite{powerlaw}) focusing on a modified Van der Pol (VdP) oscillator model, researchers observed that such transitions lead to closed trajectories resembling Hamiltonian systems like the harmonic oscillator. They analytically established this behavior for weakly nonlinear systems using multi-scale perturbative analysis. The study revealed a power law decay characteristic of center-like or quasi-periodic oscillations, with an exponent determined to be $\frac{1}{2}$.  The question was raised initially a decade ago by Bhattacharjee et al.~\cite{len3.5} through RG as a probe, but it has received less attention in the literature over the years.  Despite its mathematical nature, such studies are crucial as they can pave the way for understanding complex real-world systems, as exemplified by applications like in synthetic biology~\cite{JeffHasty,2018_RNA_Velocity_Nature,2021_RNA_Velocity_MBS}, where mathematical concepts lead to profound insights. {\color{black}However, a natural question arises regarding systems with greater nonlinearity in the damping force function, as seen in the design of birhythmic and trirhythmic oscillations. Specifically, whether the $\frac{1}{2}$ exponent remains valid under such conditions is unclear. This challenge stems from the complexity of solving amplitude dynamics equations with higher-order nonlinearities in the rate of change. To address this question, this report aims to investigate through numerical optimization algorithms. By simulating prototype systems exhibiting varying levels of nonlinearity in damping forces, the study intends to determine if and how the $\frac{1}{2}$ exponent persists in scenarios where systems suddenly transition to isochronous oscillations. This investigation is essential for advancing our understanding of rhythmic dynamics in real-world applications where such transitions can occur unexpectedly.}

Section~\ref{PD} introduces the problem formulation, outlining the mathematical perspectives as a tool of such physical questions briefly and referencing the specific case discussed in the article~\cite{powerlaw}. In Section~\ref{PD2}, multiscale perturbative method (KB) is used to convert the specific question into the mathematical frame. {\color{black}The main findings are presented using a numerical approach that extends the analysis to higher-order nonlinearities is presented in Section~\ref{RD}. This is accomplished through an optimization algorithm focusing on identifying and characterizing slowly decaying center-type oscillations.}Section~\ref{Conclusions} concludes the study, highlighting key insights and proposing future directions for research. 

\section{Description of the Problem}
\label{PD}

Considering here a generalized form of the VdP oscillator~\cite{limiso,powerlaw} instead of the traditional form where $a=1$. The parameter $a$ is crucial as it sets the amplitude for the limit cycle oscillation; for example, the general amplitude of the VdP oscillator is 2 if $a$ is set to unity. This generalized system can produce both a limit cycle and a slowly decaying center-type solution by varying the parameter $a$, unlike the traditional VdP oscillator, which only produces a limit cycle solution for nonzero values of the nonlinearity control parameter $\epsilon$. The basic equations for unit frequency is given by:$$\ddot{x}+ \epsilon (x^2-a^2)\dot{x}+ x=0;~0<\epsilon \ll 1$$ where $x$ is the time-dependent variable and the only fixed point is centered around the origin. It has been shown in Ref.~\cite{powerlaw} that there exists a decaying center-type solution for the null value of $a$ that follows a power law decay with an exponent of $\frac{1}{2}$.  Additionally, {\color{black}it remains unexplored whether this type of decaying center-type solution will follow the same power law with the same exponent in other systems, or if it is system-specific.  If it is not system-specific, is there a generalized rule for such {\color{black} class of higher order rhythmic} systems, and what will be the exponent?}

This is quite unrealistic from the limit cycle oscillation point of view, as once the trajectory falls into an isolated trajectory, it will remain forever. However, sometimes a sudden transition can occur in a system due to external perturbation, which might change the system dynamics~\cite{powerlaw,len3.5}. {\color{black}To explore this question, let us reformulate it into a mathematical modeling aspect.} A recent study~\cite{SahaCNSNS} explored the systematic design of mathematical models of the prototype VdP oscillator, described by the equation:$$\ddot{x}+ \epsilon (-a^2+x^2-\alpha x^4+\beta x^6-\gamma x^8+\delta x^{10})\dot{x}+ x=0;~0<\epsilon \ll 1$$where $\alpha, \beta, \gamma, \delta$ are system parameters, and suitable values can produce birhythmic as well as trirhythmic oscillations in families of VdP oscillators.  This prototype model is receiving significant attention in bifurcation theory, whether under the influence of noise or in noise-free conditions~\cite{Citation1,Citation2,Citation3}. It has a nice correspondence with the VdP model when all $\alpha, \beta, \gamma, \delta$ are null and $a$ is unity, showing monorhythmic behavior. Further extension shows birhythmic oscillation for nonzero values of $\alpha, \beta$, which is also known as the Kaiser model~\cite{kaiser91,k-dsr,k-y2007a}.  {\color{black}Now, the question is, for the same consideration as shown in Ref.~\cite{powerlaw}, whether there exists a decaying center-type solution for the null value of $a$ that follows a power law decay with an exponent of $\frac{1}{2}$, or whether the degree of the exponent is determined by the degree of the nonlinearities present in the damping force function. If it does not follow the exponent $\frac{1}{2}$, what will it be? Is there a general exponent for such systems that the power law will follow? The question was first raised a decade ago by Bhattacharjee et al.~\cite{len3.5} through the context of the center or limit cycle and RG as a probe, but has since received little attention in the literature. The answer will be useful for understanding the transition rate of a center-like isochronous oscillation.}

\section{Multi-scale Perturbative Method: a Tool to Reformulate a Question into a Mathematical Framework}
\label{PD2}

To answer the above questions, we will briefly describe a general multi-scale perturbative method known as the KB method~\cite{slross,len0,SahaCNSNS,lccounting}. This method allows us to convert the physical question into a mathematical framework, enabling us to visualize the problem through the amplitude dynamics of the oscillation, including the difficulties in providing a handy analytical solution. 

To begin with, let us consider the form of a generalized autonomous weakly nonlinear oscillator (often called as LLS oscillator~\cite{remickens,len4,limiso,SahaCNSNS,lccounting}) as follows:
\begin{eqnarray}
\ddot{x}+\epsilon h(x,\dot{x})+\omega^2 x=0, ~0<\epsilon\ll1,
\label{eq:KB_General}
\end{eqnarray}
where, $h(x,\dot{x})$ is a function of the position variable, $x$ and the velocity variable, $\dot{x}$ which contains nonlinearity coming from nonlinear damping or restoring force. The nonlinearity control parameter, $\epsilon$, must be within $0$ and $1$. Now, let us rewrite Eq.~\eqref{eq:KB_General} in the form as,
\begin{eqnarray}
\dot{x}&=&y, \nonumber \\
\dot{y}&=&-\omega^2 x -\epsilon h(x,\dot{x}).
\label{eq:KB_General_rewrite}
\end{eqnarray}
For $\epsilon = 0$, the above system reduces to Simple Harmonic Oscillator (SHO), with natural frequency $\omega$ having a solution, $x=r \cos(\omega t+\phi)$ and $y=-\omega~r \sin(\omega t+\phi)$, with constant amplitude, $r=\sqrt{x^2 +\frac{y^2}{\omega^2}}$ and phase, $\phi=-\omega t+ \tan^{-1}(- \frac{y}{\omega x})$ having the circular orbit of period $\frac{2 \pi}{\omega}$.

When $\epsilon \neq 0$, indicating the activation of nonlinear terms and a change in the vector field by an amount $\epsilon$, the orbits of the system become nearly circular for small $\epsilon$.  These orbits approximately repeat every $\frac{2\pi}{\omega}$, where $\omega$ denotes a characteristic frequency.  Even small variations in $\epsilon$ can noticeably alter the period of these orbits.  This raises the non-trivial question of how to characterize the amplitude-phase dynamics under these conditions.

The approach involves combining power consumption law with averaging theory~\cite{strogatz,len0}.  In the context of a SHO where total energy is conserved, the system's energy remains unchanged over one cycle.  Introducing $\epsilon$ results in occasional energy input or damping, adding an intriguing dynamic aspect to the system.  These fluctuations can lead to the formation of limit cycles, influenced by the interplay between these effects.  By explicitly considering the amplitude $r$ and phase $\phi$ as time-dependent quantities, we gain insights into the evolution of solutions under weak nonlinearity.

Therefore, considering $x \approx r(t) \cos(\omega t+\phi(t))$ and $y \approx -\omega~r(t) \sin(\omega t+\phi(t))$  with $r(t) \approx \sqrt{x^2+\frac{y^2}{\omega^2}}$ and $\phi(t) \approx -\omega t+ \tan^{-1} \left(-\frac{y}{\omega x}\right)$, then one can have the rate of change of amplitude and phase variables with respect to time as,
\begin{eqnarray}
\dot{r}(t)&=&\frac{\epsilon~h}{\omega} \sin(\omega t+\phi(t)),\nonumber\\
\dot{\phi}(t)&=&\frac{\epsilon~h}{\omega~r(t)} \cos(\omega t+\phi(t)).
\end{eqnarray}
This implies that the time derivatives of amplitude and phase are of $O(\epsilon)$. Now, if $\overline{U}(t)$ be a running  average of a time dependent function $U$ defined as,
\begin{eqnarray}
\overline{U}(t)= \frac{\omega}{2 \pi} \int_{0}^{\frac{2 \pi}{\omega}} U(s)~ds,
\end{eqnarray} 
then from the fundamental theorem of calculus it is observed that $\dot{\overline{U}}=\overline{\dot{U}}$. By applying this averaging trick to the time dependent amplitude and phase equations for each cycle (as the trial solution is taken approximately periodic) then we have,
\begin{eqnarray}
\dot{\overline{r}} &= \langle \frac{\epsilon ~ h(x,y)}{\omega} ~ \sin (\omega t+\phi(t)) \rangle_t,\nonumber\\
\dot{\overline{\phi}} &= \langle \frac{\epsilon ~ h(x,y)~}{\omega~r(t)} \cos (\omega t+\phi(t))\rangle_t.
\end{eqnarray}
As $\dot{r}(t)$ and $\dot{\phi}(t)$ are of $O(\epsilon)$ then we may set the perturbation on $r$ and $\phi$ over each cycle as,
\begin{eqnarray}
r(t) &=&\overline{r}+O(\epsilon), \nonumber\\
\phi(t) &=&\overline{\phi}+O(\epsilon),
\end{eqnarray}
where, $\overline{r}$ and $\overline{\phi}$ are very weakly time dependent so that the error can be negligible. Finally, from the above consideration, one can obtain,
\begin{eqnarray}
\dot{\overline{r}} &=&\langle \frac{\epsilon ~ h(\overline{r} \cos(\omega t+\overline{\phi}),-\omega~\overline{r} \sin(\omega t+\overline{\phi})}{\omega} ~ \sin (\omega t+\overline{\phi}) \rangle_t \nonumber\\
&=& \varphi_1 (\overline{r},\overline{\phi}) \nonumber\\
\dot{\overline{\phi}} &=& \langle \frac{\epsilon ~ h(\overline{r} \cos(\omega t+\overline{\phi}),-\omega~\overline{r} \sin(\omega t+\overline{\phi})}{\omega~\overline{r}} \cos (\omega t+\overline{\phi})\rangle_t \nonumber\\
&=& \varphi_2 (\overline{r},\overline{\phi})
\label{eq:KB_General_rewrite_AmpPhEq}
\end{eqnarray}
where, $O(\epsilon^2)$ terms can be neglected as first order approximation. The above systems will not be in a coupled form for autonomous set of equation in general, but for the non-autonomous case or in the case of more than one potential well we may have coupling between amplitude and phase flow variables in the flow equation and that kind of situation may provide very complex phenomena which is quite harder to solve as well as to get back the original system of equations which can be fixed by considering further approximations~\cite{physicad}.

In the following discussion, if we examine the amplitude dynamics of a slightly modified VdP oscillator, we get $\dot{\overline{r}} = - \frac{\epsilon \overline{r}}{8} \left( \overline{r}^2-4 a^2 \right)$. The phase equation remains constant due to the autonomous structure of the system and its single well potential. If a sudden change occurs, setting $a$ to zero, the amplitude dynamics become $\dot{\overline{r}} = - \frac{\epsilon \overline{r}^3}{8}$. This equation is solvable and yields the exact solution $\overline{r}=\frac{2}{\sqrt{\frac{4}{\overline{r}_0^2}+\epsilon  t}}$. This solution can be approximated as $\overline{r} \propto t^{-\frac{1}{2}}$.

This allows us to identify the nature and rule of the exponent that the power law follows for a decaying center-type solution, which is visualizable through amplitude dynamics. {\color{black}The question then arises: what happens in the generalized birhythmic and trirhythmic model when such a sudden change sets $a$ to zero? In this case, the amplitude equation reads $\dot{\overline{r}}=\frac{\epsilon  \overline{r}}{1024} \left(64 \alpha  \overline{r}^4-40 \beta  \overline{r}^6+28 \gamma  \overline{r}^8-21 \delta  \overline{r}^{10}-128 \overline{r}^2\right)$. This equation reduces to the previous case when $\alpha, \beta, \gamma, \delta$ are all zero. } 

{\color{black}In addition to the VdP class, we consider an alternative class called the Rayleigh class~\cite{SahaCNSNS}, whose basic equation reads: $$\ddot{x}+ \epsilon (-a^2+\dot{x}^2-\alpha \dot{x}^4+\beta \dot{x}^6-\gamma \dot{x}^8+\delta \dot{x}^{10})\dot{x}+ x=0;~0<\epsilon \ll 1.$$ 

Upon applying the KB perturbative approach and carefully accounting for an external nuisance, we derive the amplitude equation:  $\dot{\overline{r}}=\frac{\epsilon  \overline{r}}{1024} \left(320 \alpha  \overline{r}^4-280 \beta  \overline{r}^6+252 \gamma  \overline{r}^8-231 \delta  \overline{r}^{10}-384 \overline{r}^2\right)$.}

{\color{black}The main difficulty is the lack of an analytical method to solve this type of nonlinear ordinary differential equation, which prevents us from determining {\color{black} an approximate} exponent in terms of time, even when only $\alpha$ and $\beta$ are nonzero, resulting in birhythmicity. To address this problem, we have adopted a numerical approach, which is discussed below.}

\section{Algorithm and Corresponding Numerical Results}
\label{RD}
To address this issue, we adopted several numerical algorithms sequentially using open-source Python libraries. {\color{black}Also, for simplification,  we set $\gamma=\delta=0$.} {\color{black}First, we used `scipy.integrate' to solve the complex, highly nonlinear ordinary differential equation. Then, we used `scipy.optimize' to fit a curve of the form $\propto t^{-\frac{1}{n}}$ to the numerical solutions obtained for thousand of random initial values of $r_0$ within the interval $(2.4, 4.91)$. This interval was chosen because it includes the amplitudes of the {\color{black}bi}rhythmic oscillator, which are $2.63902$, $3.96164$, and $4.83953$~\cite{SahaCNSNS}. These amplitudes represent possible starting points for decaying center-type oscillations in case of sudden transitions. }{\color{black} A similar consideration for $r_0$ is applied to the Rayleigh class of oscillations within the interval $(1.5,2.6)$,  as the birhythmic oscillation admits amplitudes of $1.69091$, $2.03334$, and $2.51274$~\cite{SahaCNSNS}.}

Finally, we used `sklearn.metrics' to find the least mean squared error among the best-fitted curves for various initial conditions. For the birhythmic {\color{black}VdP class}, $\alpha$ and $\beta$ were set to $0.144$ and $0.005$, respectively~\cite{k-dsr,SahaCNSNS} {\color{black} and for the birhythmic Rayleigh class, $\alpha$ and $\beta$ were set to $0.0.285272$ and $0.0244993$, respectively~\cite{SahaCNSNS}. } {\color{black} In all cases, we observed that the decaying center-type oscillation followed a power law with an exponent of $\frac{1}{3}$, indicating that $n=3$ provided the best fit, with the mean squared error being minimal for various random initial conditions within the interval. }

\begin{figure}[htbp]
\centering
\includegraphics[width=0.5\textwidth]{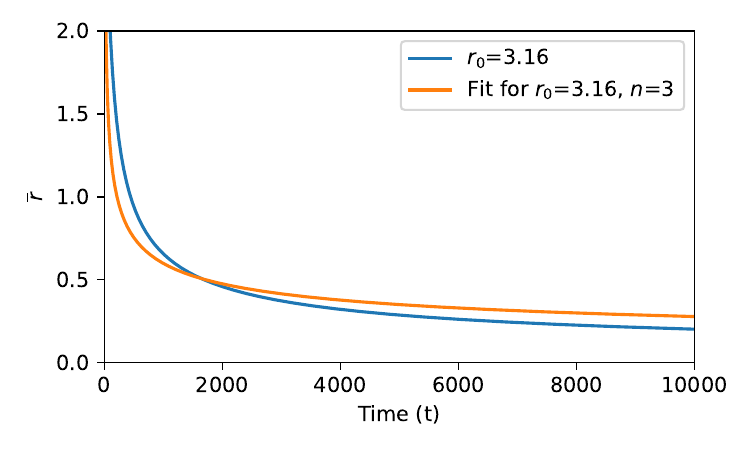}\\
\includegraphics[width=0.5\textwidth]{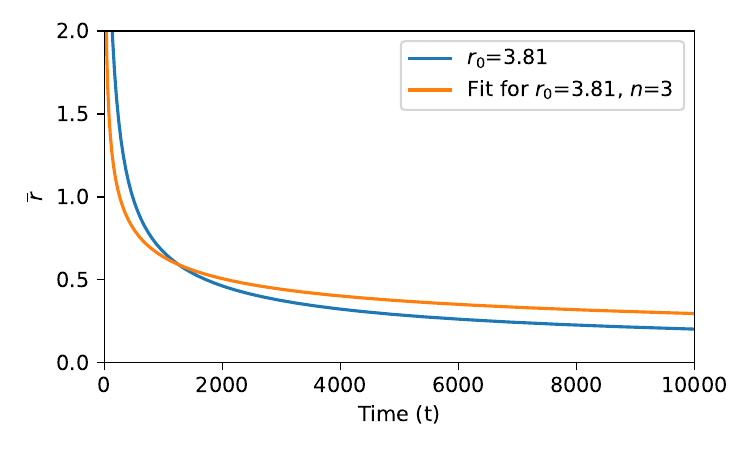}
\caption{Numerical solution and fitted functions for different initial conditions for the Birhythmic case (VdP class) $(\alpha=0.144,~\beta=0.005)$: $r_0=$3.16 (top) and $r_0=$3.81 (bottom).}
 \label{fig:Birhythmic}
\end{figure}

\begin{figure}[htbp]
\centering
\includegraphics[width=0.5\textwidth]{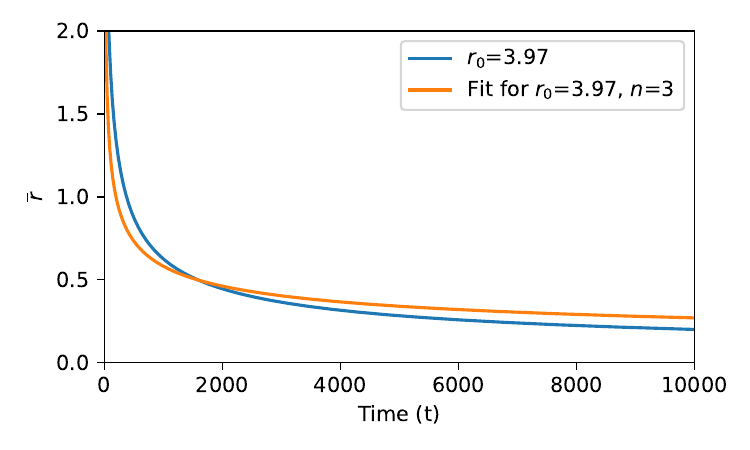}\\
\includegraphics[width=0.5\textwidth]{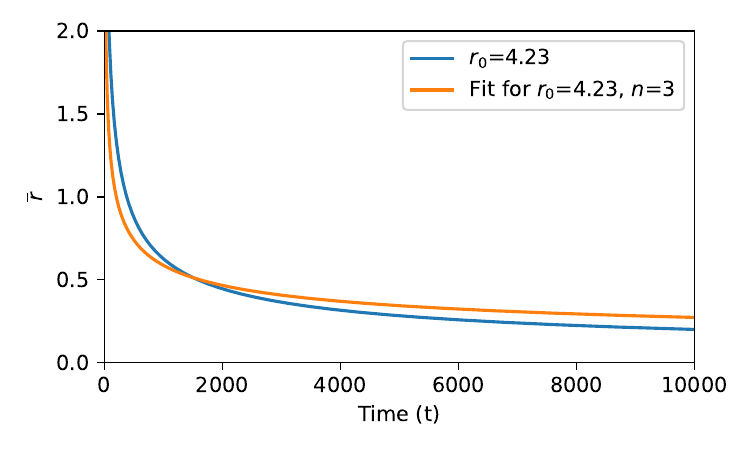}
\caption{Numerical solution and fitted functions for different initial conditions for the Monorhythmic case (VdP class) $(\alpha=\beta=\gamma=\delta=0)$: $r_0=$3.97 (top) and $r_0=$4.23 (bottom).}
 \label{fig:Monorhythmic}
\end{figure}

\begin{figure}[htbp]
\centering
\includegraphics[width=0.5\textwidth]{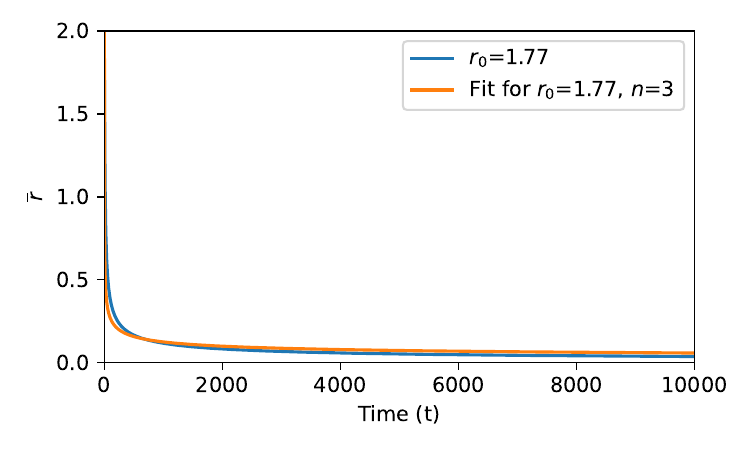}\\
\includegraphics[width=0.5\textwidth]{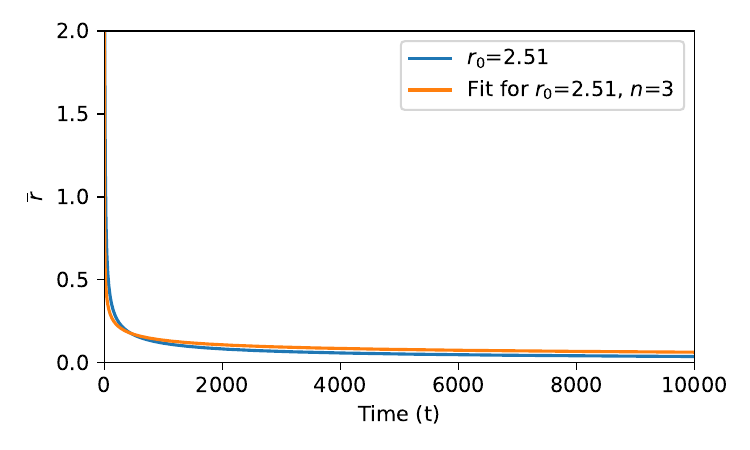}
\caption{\color{black}Numerical solution and fitted functions for different initial conditions for the Birhythmic case (Rayleigh class) $(\alpha=0.285272,~\beta=0.0244993)$: $r_0=$1.77 (top) and $r_0=$2.51 (bottom).}
 \label{fig:Birhythmic_2}
\end{figure}

\begin{figure}[htbp]
\centering
\includegraphics[width=0.5\textwidth]{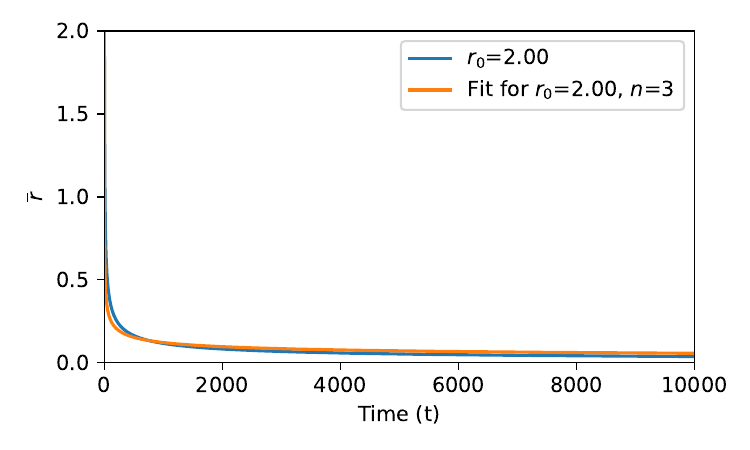}\\
\includegraphics[width=0.5\textwidth]{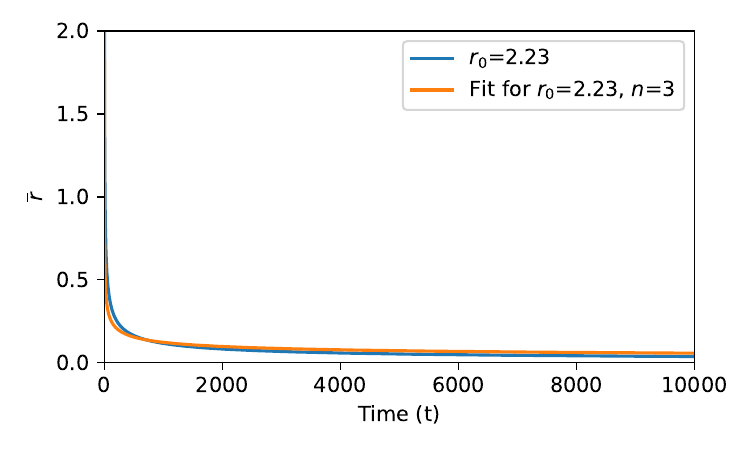}
\caption{\color{black}Numerical solution and fitted functions for different initial conditions for the Monorhythmic case (Rayleigh class) $(\alpha=\beta=\gamma=\delta=0)$: $r_0=$2.0 (top) and $r_0=$2.23 (bottom).}
 \label{fig:Monorhythmic_2}
\end{figure}

The figures labeled Fig.~\ref{fig:Birhythmic} and Fig.~{\color{black}\ref{fig:Birhythmic_2}} provide numerical solutions and fitted functions for various random initial conditions in both birhythmic cases. 
Interestingly, in both scenarios, the optimal fitted index, denoted as $n$ and corresponding to the power law exponent of $\frac{1}{3}$, is observed to be 3. This consistency extends even to the monorhythmic case (Fig.~\ref{fig:Monorhythmic} and {\color{black}\ref{fig:Monorhythmic_2}}), akin to a slightly VdP or Rayleigh type oscillator, where the fitted index is also 3. This alignment is noteworthy as previous studies relied on an approximation of the original analytical solution, which might not always hold true when compared against precise numerical solutions. Consequently, for generalized nonlinear systems with prominently nonlinear terms in the damping force function, the power law emerges as a standard method for modeling the decay of center-type solutions following sudden changes, with the typical exponent being $\frac{1}{3}$.


\section{Conclusions and Future Directions}
\label{Conclusions}
In conclusion, by suitably adopting the KB averaging method in multi-scale perturbation theory, we have established a general rule for the decaying center-like power law behavior, characterized by an exponent of $\frac{1}{3}$ {\color{black} (at least for the birhythmic variations of the classes)}. The investigation began with a physical question regarding the higher-order nonlinearity involved in a damping force function, which admits birhythmic as well as trirhythmic behavior under a sudden transition to a center-type solution. Due to the higher-order nonlinearity in dissipation, the solution is not as straightforward as in the Hamiltonian case.

We successfully translated the physical problem into a mathematical framework by employing the KB averaging approach to derive the amplitude-phase dynamics. Using numerical algorithms from open-source Python libraries, we optimized the exponent by fitting curves for thousands of random initial conditions. For all these conditions, we observed a consistent general rule of decaying power law behavior with an exponent of $\frac{1}{3}$, regardless of whether the system initially exhibited monorhythmic {\color{black}or} birhythmic behavior.

{\color{black}A limitation of this approach is the careful selection of random initial conditions close to the rhythmic closed orbits. If the transition occurs, the decaying would start near such cycles for non-zero values of time; otherwise, there is a risk of divergent behavior. As known, the exact center solution for Hamiltonian systems is sensitive to initial conditions, as the initial values determine the amplitude of motion. Therefore, this method partially mirrors the nature of an exact center solution, as the condition for having a center or decaying-center-like system is the same--namely, the zero value of the constant damping force function.}

As a future direction, {\color{black}investigating the behavior of the described systems with $\gamma ~\text{and}~ \delta$ could be valuable in formulating a generalized rule for trirhythmic cases.} It would also be quite interesting to investigate the behavior of the described systems in the presence of a periodic external force. This has practical applications in achieving self-sustained oscillations, particularly when the systems exhibit exponential decay and have a stable fixed point~\cite{dsrrayleigh}. Additionally, examining the behavior under periodic modulation of the damping force function~\cite{dsrrayleigh,physicad} during such transitions would be valuable. It would be intriguing to determine whether resonance competes with the slowly decaying center-like or quasi-periodic solution and the periodic modulation, and if so, whether this resonance could help control such types of power law decay. This study would further elucidate the system's future dynamics under such modulation. {We believe this key result would be beneficial for understanding the multi-rhythmic nature of biological systems and engineering systems (through nonlinear circuits).


\section*{Acknowledgement}
\noindent
{SS acknowledges Prof. G. Gangopadhyay (SNBNCBS, Kolkata), Prof. D. S. Ray (IACS Kolkata), Prof. D. Mondal (IIT Tirupati), Prof. P. Ghosh (IISER TVM),Prof. S. Kar (IIT Bombay), Dr. S. Ray (IISER BBSR), Dr. S. Chakraborty (FAU Germany), Prof. S. Chakraborty (IIT Kanpur), Prof. J. N. Tegner (KAUST, KSA) and Dr. A. Pandey (GU, Noida) for their valuable support and insightful discussions.  SS would also like to thank ChatGPT for assisting with language improvements and editing suggestions.} 

\bibliographystyle{unsrt}
\bibliography{References-ssaha}{}

\end{document}